# A characterization of trace zero symmetric nonnegative 5x5 matrices

Oren Spector

## Abstract

The problem of determining necessary and sufficient conditions for a set of real numbers to be the eigenvalues of a symmetric nonnegative matrix is called the symmetric nonnegative inverse eigenvalue problem (SNIEP). In this paper we solve SNIEP in the case of trace zero symmetric nonnegative $5 \times 5$ matrices.

## 1. Introduction

The problem of determining necessary and sufficient conditions for a set of complex numbers to be the eigenvalues of a nonnegative matrix is called the nonnegative inverse eigenvalue problem (NIEP). The problem of determining necessary and sufficient conditions for a set of real numbers to be the eigenvalues of a nonnegative matrix is called the real nonnegative inverse eigenvalue problem (RNIEP). The problem of determining necessary and sufficient conditions for a set of real numbers to be the eigenvalues of a symmetric nonnegative matrix is called the symmetric nonnegative inverse eigenvalue problem (SNIEP). All three problems are currently unsolved in the general case.

Loewy and London [6] have solved NIEP in the case of $3 \times 3$ matrices and RNIEP in the case of $4 \times 4$ matrices. Moreover, RNIEP and SNIEP are the same in the case of $n \times n$ matrices for $n \leq 4$. This can be seen from papers by Fielder [2] and Loewy and London [6]. However, it has been shown by Johnson et al. [3] that RNIEP and SNIEP are different in general. More results about the general NIEP, RNIEP and SNIEP can be found in [1]. Other results about SNIEP in the case $n = 5$ can be found in [7] and [8].

Reams [9] has solved NIEP in the case of trace zero nonnegative $4 \times 4$ matrices. Laffey and Meehan [4] have solved NIEP in the case of trace zero nonnegative $5 \times 5$ matrices. In this paper we solve SNIEP in the case of trace zero symmetric nonnegative $5 \times 5$ matrices.

The paper is organized as follows: In Section 2 we present some notations. We also give some basic necessary conditions for a spectrum to be realized by a trace zero symmetric nonnegative $5 \times 5$ matrix. In Section 3 we state our main result without proof. In Section 4 we present some preliminary results that are needed for the proof of the main results. Finally, in Section 5 we prove our main result.

## 2. Notations and some necessary conditions

Let $\mathfrak{R}$ be the set of trace zero nonnegative $5 \times 5$ matrices with real eigenvalues. Let $A \in \mathfrak{R}$. We call a spectrum $\sigma = \sigma(A)$ a *normalized spectrum* if $\sigma = (\lambda_1, \lambda_2, \lambda_3, \lambda_4, \lambda_5)$ with $1 = \lambda_1 \geq \lambda_2 \geq \lambda_3 \geq \lambda_4 \geq \lambda_5 \geq -1$. Since $A$ has a zero trace $\sum_{i=1}^{5} \lambda_i = 0$.

Let $\widetilde{\mathfrak{R}}$ be the set of trace zero symmetric nonnegative $5 \times 5$ matrices. Then $\widetilde{\mathfrak{R}} \subset \mathfrak{R}$.

Notice that for a normalized spectrum $\sigma = (\lambda_1, \lambda_2, \lambda_3, \lambda_4, \lambda_5)$:

$$\lambda_5 = \frac{1}{4}(\lambda_5 + \lambda_5 + \lambda_5 + \lambda_5) \leq \frac{1}{4}(\lambda_2 + \lambda_3 + \lambda_4 + \lambda_5) = -\frac{1}{4}\lambda_1 = -\frac{1}{4},$$

$$\lambda_2 + \lambda_3 + \lambda_4 = -\lambda_1 - \lambda_5 = -1 - \lambda_5 \geq -1 + \frac{1}{4} = -\frac{3}{4},$$

$$\lambda_2 + \lambda_3 + \lambda_4 = -\lambda_1 - \lambda_5 \leq -1 + 1 = 0.$$

Let $d = \lambda_2 + \lambda_3 + \lambda_4$. Then $d \in \left[-\frac{3}{4}, 0\right]$, $\lambda_4 = d - \lambda_2 - \lambda_3$ and $\lambda_5 = -d - 1$.

For the rest of this paper we shall use $x$ instead of $\lambda_2$ and $y$ instead of $\lambda_3$. We shall deal mainly with spectra of the form $\sigma = (1, x, y, d - x - y, -d - 1)$.

We are interested in finding necessary and sufficient conditions for $\sigma = (1, x, y, d - x - y, -d - 1)$ to be a normalized spectrum of a matrix $A \in \widetilde{\mathfrak{R}}$.

We start by looking at the necessary conditions that a normalized spectrum $\sigma = (1, x, y, d - x - y, -d - 1)$ imposes.

From $1 \geq x \geq y \geq d - x - y \geq -d - 1$ we get:

$$x \leq 1,$$

$$y \leq x,$$

$$y \leq -x + 2d + 1,$$

$$y \geq \frac{1}{2}(d - x).$$

Let the following points on the $xy$ plain be defined:

$$A = \left(\frac{1}{3}d, \frac{1}{3}d\right),$$

$$B = \left(d + \frac{1}{2}, d + \frac{1}{2}\right),$$

$$C = (3d + 2, -d - 1),$$

$$D = (1,2d),$$

$$E = \left(1, \frac{1}{2}d - \frac{1}{2}\right),$$

$$F = (d + 1, d),$$

$$G = \left(d + 1, -\frac{1}{2}\right),$$

$$H = \big(f(d), f(d)\big),$$

$$I = \left(d + \frac{1}{2} + g(d), d + \frac{1}{2} - g(d)\right),$$

$$J = (2d + 1,0),$$

$$O = (0,0),$$

where,

$$f(d) = \frac{2}{3}d - \frac{\sqrt[3]{4}d^2}{3\sqrt[3]{r(d)}} - \frac{\sqrt[3]{r(d)}}{3\sqrt[3]{4}},$$

$$r(d) = 4d^3 + 27d^2 + 27d + 3\sqrt{3}\sqrt{d^2(d + 1)(8d^2 + 27d + 27)},$$

$$g(d) = \sqrt{\frac{\left(d - \left(-\frac{3}{4} + \frac{\sqrt{5}}{4}\right)\right)\left(d - \left(-\frac{3}{4} - \frac{\sqrt{5}}{4}\right)\right)}{2d + 1}}.$$

A few notes are in order:

1. $r(d)$ is real for $d \geq -1$ so $f(d)$ is real for $d \geq -1$ (except perhaps when $d = 0$),
2. $g(d)$ is real for $d \geq -\frac{3}{4} + \frac{\sqrt{5}}{4}$.

For $d \in \left[-\frac{3}{4}, -\frac{1}{3}\right]$ the above normalization conditions form the triangle $ABC$ in the $xy$ plain. Note that for $d = -\frac{3}{4}$ the triangle becomes a single point with coordinates $\left(-\frac{1}{4}, -\frac{1}{4}\right)$.

For $d \in \left[-\frac{1}{3}, 0\right]$ the conditions form the quadrangle $ABDE$ in the $xy$ plain.

Another necessary condition we shall use is due to McDonald and Neumann [8] (Lemma 4.1).

**Theorem 1 (MN)**: Let $A$ be a $5 \times 5$ irreducible nonnegative symmetric matrix with a spectrum $\sigma(A) = (\lambda_1, \lambda_2, \lambda_3, \lambda_4, \lambda_5)$ such that $\lambda_1 \geq \lambda_2 \geq \lambda_3 \geq \lambda_4 \geq \lambda_5$. Then $\lambda_2 + \lambda_5 \leq trace(A)$.

Loewy and McDonald [7] extended this result to any $5 \times 5$ nonnegative symmetric matrix (not just irreducible).

In our case $\sigma = (1, x, y, d - x - y, -d - 1)$ so we get the necessary condition $x \leq d + 1$.

The MN condition is already met when $d \in \left[-\frac{3}{4}, -\frac{1}{2}\right]$, while it limits the realizable set to the quadrangle $ABFG$ when $d \in \left[-\frac{1}{2}, 0\right]$.

For $d \in \left[-\frac{3}{4} + \frac{\sqrt{5}}{4}, 0\right]$ let $P$ be the 5-vertex shape $AHIFG$ such that all its edges but $HI$ are straight lines and the edge $HI$ is described by the curve $\gamma: [0,1] \to \mathbb{R}^2$, where

$$\gamma(t) = \big(x(t), h(x(t))\big),$$

$$x(t) = (1 - t)f(d) + t\left(d + \frac{1}{2} + g(d)\right),$$

$$h(t) = -\frac{1}{2}(t - d) + \frac{1}{2}\sqrt{\frac{t^3 + dt^2 - d^2t - 4d - 4d^2 - d^3}{t - d}}.$$

We shall show later that $h(t)$ is real for $t \in \left[f(d), d + \frac{1}{2} + g(d)\right]$, $h(f(d)) = f(d)$ and $h\left(d + \frac{1}{2} + g(d)\right) = d + \frac{1}{2} - g(d)$, so that $\gamma(0) = H, \gamma(1) = I$ which makes $P$ well defined.

Let

$$\bar{r}(d) = 4d^3 + 27d^2 + 27d - 3\sqrt{3}\sqrt{d^2(d+1)(8d^2+27d+27)}.$$

For $d < 0$

$$\bar{r}(d) = d\left(4d^2 + 27d + 27 + 3\sqrt{3}\sqrt{(d+1)(8d^2+27d+27)}\right),$$

$$r(d) = \frac{r(d)\bar{r}(d)}{\bar{r}(d)} = \frac{16d^5}{4d^2 + 27d + 27 + 3\sqrt{3}\sqrt{(d+1)(8d^2+27d+27)}}.$$

Note that $\lim_{d\to 0^-} \frac{r(d)}{d^5} = \frac{16}{54}$ which means $r(d) = O(d^5)$, so $\lim_{d\to 0^-} f(d) = 0$. Moreover, $g(0) = \frac{1}{2}$ and when $d = 0$ we have $h(t) = 0$ for $t \geq 0$. Therefore, when $d = 0$ we get $H = O = A$ and $I = J = F = D$ so $P$ becomes the triangle $AFG = OJG$.

See Appendix A for the orientation of the points $A$–$J$,$O$ in the plain for different values of $d$.

## 3. Statement of main result

The following theorem completely solves SNIEP in the case of trace zero symmetric nonnegative $5 \times 5$ matrices.

**Theorem 2 (Main result)**: A necessary and sufficient condition for $\sigma = (1, x, y, d - x - y, -d - 1)$ to be a normalized spectrum of a matrix $A \in \widetilde{\mathfrak{R}}$ is:

a) $(x, y)$ lies within the triangle $ABC$ for $d \in \left[-\frac{3}{4}, -\frac{1}{2}\right]$,
b) $(x, y)$ lies within the quadrangle $ABFG$ for $d \in \left[-\frac{1}{2}, -\frac{3}{4} + \frac{\sqrt{5}}{4}\right]$,
c) $(x, y)$ lies within $P$ for $d \in \left[-\frac{3}{4} + \frac{\sqrt{5}}{4}, 0\right]$.

An immediate corollary of Theorem 2 is:

**Theorem 3**: Let $\sigma = (\lambda_1, \lambda_2, \lambda_3, \lambda_4, \lambda_5)$ and let $s_k = \sum_{i=1}^{5} {\lambda_i}^k$. Suppose $\lambda_1 \geq \lambda_2 \geq \lambda_3 \geq \lambda_4 \geq \lambda_5$. Then $\sigma$ is a spectrum of a matrix $A \in \widetilde{\mathfrak{R}}$ if and only if the following conditions hold:

1. $s_1 = 0$,
2. $s_3 \geq 0$,
3. $\lambda_2 + \lambda_5 \leq 0$.

## 4. Preliminary results

In order to prove Theorem 2 we need several results.

The first result is due to Fiedler [2], which extended a result due to Suleimanova [10]:

**Theorem 4 (Fiedler)**: Let $\lambda_1 \geq 0 \geq \lambda_2 \geq \cdots \geq \lambda_n$ and $\sum_{i=1}^{n} \lambda_i \geq 0$. Then there exists a symmetric nonnegative $n \times n$ matrix $A$ with a spectrum $\sigma = (\lambda_1, \lambda_2, \dots, \lambda_n)$, where $\lambda_1$ is its Perron eigenvalue.

The second result is due to Loewy [5]. Since this result is unpublished we shall give Loewy's proof.

**Theorem 5 (Loewy)**: Let $n \geq 4$, $\lambda_1 \geq \lambda_2 \geq 0 \geq \lambda_3 \geq \cdots \geq \lambda_n$, $\sum_{i=1}^{n} \lambda_i \geq 0$. Suppose $K_1, K_2$ is a partition of $\{3,4, \dots, n\}$ such that $\lambda_1 \geq -\sum_{i \in K_1} \lambda_i \geq -\sum_{i \in K_2} \lambda_i$. Then there exists a nonnegative symmetric $n \times n$ matrix $A$ with a spectrum $\sigma = (\lambda_1, \lambda_2, \dots, \lambda_n)$, where $\lambda_1$ is its Perron eigenvalue.

For the proof of Theorem 5 we shall need another result of Fiedler [2]:

**Theorem 6 (Fiedler)**: Let the following conditions hold:

1. $(\alpha_1, \alpha_2, \dots, \alpha_m)$ is realizable by a nonnegative symmetric $m \times m$ matrix with $\alpha_1$ as its Perron eigenvalue,
2. $(\beta_1, \beta_2, \dots, \beta_n)$ is realizable by a nonnegative symmetric $n \times n$ matrix with $\beta_1$ as its Perron eigenvalue,
3. $\alpha_1 \geq \beta_1$,
4. $\varepsilon \geq 0$.

Then $(\alpha_1 + \varepsilon, \beta_1 - \varepsilon, \alpha_2, \dots, \alpha_m, \beta_2, \dots, \beta_n)$ is realizable by a nonnegative symmetric $(m+n) \times (m+n)$ matrix with $\alpha_1 + \varepsilon$ as its Perron eigenvalue.

**Proof of Theorem 5 (Loewy)**: Let $\varepsilon = \lambda_1 + \sum_{i \in K_1} \lambda_i$. Then by the assumptions $\varepsilon \geq 0$. We deal with two cases.

If $\lambda_1 - \varepsilon \geq \lambda_2 + \varepsilon$ then by the assumptions the set $\{\lambda_1 - \varepsilon\} \cup \{\lambda_i | i \in K_1\}$ meets the conditions of Theorem 4. Thus, there exists a symmetric nonnegative matrix whose eigenvalues are $\{\lambda_1 - \varepsilon\} \cup \{\lambda_i | i \in K_1\}$ with $\lambda_1 - \varepsilon$ as its Perron eigenvalue. Similarly, the set $\{\lambda_2 + \varepsilon\} \cup \{\lambda_i | i \in K_2\}$ meets the conditions of Theorem 4, so there exists a symmetric nonnegative matrix whose eigenvalues are $\{\lambda_2 + \varepsilon\} \cup \{\lambda_i | i \in K_2\}$, where $\lambda_2 + \varepsilon$ is its Perron eigenvalue. As $\lambda_1 - \varepsilon \geq \lambda_2 + \varepsilon$ and $\varepsilon \geq 0$ all the conditions of Theorem 6 are met and therefore the set $\{\lambda_1, \lambda_2, \dots, \lambda_n\}$ is realizable by a nonnegative symmetric $n \times n$ matrix with $\lambda_1$ as its Perron eigenvalue.

If $\lambda_1 - \varepsilon < \lambda_2 + \varepsilon$ then let $\delta = \frac{1}{2}(\lambda_1 - \lambda_2) \geq 0$. We have $\varepsilon > \frac{1}{2}(\lambda_1 - \lambda_2) = \delta$ and $\lambda_1 - \delta = \lambda_2 + \delta = \frac{1}{2}(\lambda_1 + \lambda_2) \geq 0$. Also $\lambda_1 - \delta + \sum_{i \in K_1} \lambda_i = \varepsilon - \delta > 0$, so the set $\{\lambda_1 - \delta\} \cup \{\lambda_i | i \in K_1\}$ meets the conditions of Theorem 4. We infer that there exists a symmetric nonnegative matrix whose eigenvalues are $\{\lambda_1 - \delta\} \cup \{\lambda_i | i \in K_1\}$ with $\lambda_1 - \delta$ as its Perron eigenvalue. Also,

$$\lambda_2 + \delta + \sum_{i \in K_2} \lambda_i \geq \lambda_2 + \delta + \sum_{i \in K_1} \lambda_i = \lambda_2 + \delta + \varepsilon - \lambda_1 = \delta + \varepsilon - 2\delta = \varepsilon - \delta > 0.$$

Therefore, the set $\{\lambda_2 + \delta\} \cup \{\lambda_i | i \in K_2\}$ meets the conditions of Theorem 4. Thus, there exists a symmetric nonnegative matrix whose eigenvalues are $\{\lambda_2 + \delta\} \cup \{\lambda_i | i \in K_2\}$ with $\lambda_2 + \delta$ as its Perron eigenvalue. Again, all the conditions of Theorem 6 are met and therefore the set $\{\lambda_1, \lambda_2, \dots, \lambda_n\}$ is realizable by a nonnegative symmetric $n \times n$ matrix with $\lambda_1$ as its Perron eigenvalue.

This completes the proof of the theorem.

Let $d \in \left[-\frac{1}{2}, 0\right]$ and let $\sigma = (\lambda_1, \lambda_2, \lambda_3, \lambda_4, \lambda_5)$. Define $s_k = \sum_{i=1}^{5} {\lambda_i}^k$. A necessary condition for $\sigma$ to be a spectrum of a nonnegative $5 \times 5$ matrix $A$ is $s_k \geq 0$ for $k = 1,2,3, \dots$. This easily follows from the fact that $s_k = trace(A^k)$.

We shall investigate the properties of $s_k$ for $\sigma = (1, x, y, d - x - y, -d - 1)$ when $(x, y)$ lies within the triangle $OBJ$. Note that for $d \in \left[-\frac{1}{2}, 0\right]$ the triangle $OBJ$ is contained in the triangle $ABC$.

**Lemma 1**: Let $d \in \left[-\frac{1}{2}, 0\right]$ and a positive integer $k$ be set. Let $s_k = 1 + x^k + y^k + (d - x - y)^k + (-d - 1)^k$ where $(x, y)$ lies within the triangle $OBJ$. The following table summarizes the minimum and maximum values of $s_k$ over the triangle $OBJ$:

| $k$ | Even | Odd |
|---|---|---|
| Minimum achieved at | $(0,0)$ | $\left(d + \frac{1}{2}, d + \frac{1}{2}\right)$ |
| Minimum value | $1 + d^k + (-d - 1)^k$ | $1 + 2\left(d + \frac{1}{2}\right)^k + 2(-d - 1)^k$ |
| Maximum achieved at | $(2d + 1, 0)$ | $(0,0)$ |
| Maximum value | $1 + (2d + 1)^k + 2(-d - 1)^k$ | $1 + d^k + (-d - 1)^k$ |

**Proof**: First we note that for $d = -\frac{1}{2}$ the triangle $OBJ$ becomes a single point with coordinates $(0,0)$. In this case $s_k = 1 + 2\left(-\frac{1}{2}\right)^k$, so for even $k$ we have

$$s_k = 1 + \left(\frac{1}{2}\right)^{k-1} = 1 + \left(-\frac{1}{2}\right)^k + \left(-\left(-\frac{1}{2}\right) - 1\right)^k$$
$$= 1 + (2\left(-\frac{1}{2}\right) + 1)^k + 2\left(-\left(-\frac{1}{2}\right) - 1\right)^k,$$

and for odd $k$ we have

$$s_k = 1 - \left(\frac{1}{2}\right)^{k-1} = 1 + 2\left(-\frac{1}{2} + \frac{1}{2}\right)^k + 2\left(-\left(-\frac{1}{2}\right) - 1\right)^k$$
$$= 1 + \left(-\frac{1}{2}\right)^k + \left(-\left(-\frac{1}{2}\right) - 1\right)^k.$$

The triangle $OBJ$ is characterized by the following inequalities:

$y \leq x,$

$y \leq -x + 2d + 1,$

$y \geq 0.$

We investigate $s_k$ within the triangle by looking at lines of the form $y = wx$, where $w \in [0,1]$. These lines cover the entire triangle. For a given $w$ we have $x \in \left[0, \frac{2d+1}{w+1}\right]$.

$$s_k = 1 + x^k + y^k + (d - x - y)^k + (-d - 1)^k$$
$$= 1 + (w^k + 1)x^k + (d - (w + 1)x)^k + (-d - 1)^k,$$

$$\frac{\partial s_k}{\partial x} = k(w^k + 1)x^{k-1} - k(w + 1)(d - (w + 1)x)^{k-1}.$$

Consider $d = 0$. As $x \geq 0$ and $w \geq 0$, for even $k$ we have

$$\frac{\partial s_k}{\partial x} = k\Big((w^k + 1) + (w + 1)^k\Big)x^{k-1} \geq 0,$$

and for odd $k$ we have

$$\frac{\partial s_k}{\partial x} = k\Big((w^k + 1) - (w + 1)^k\Big)x^{k-1} \leq 0.$$

Next assume $-\frac{1}{2} \leq d < 0$. The derivative is zero when

$$k(w^k + 1)x^{k-1} = k(w + 1)(d - (w + 1)x)^{k-1}.$$

We know $d - (w + 1)x \neq 0$. Otherwise, we must have $x = 0$ and therefore $d = 0$, which is a contradiction to our assumption. Then, as $0 \leq w \leq 1$, we have

$$\left(\frac{x}{d - (w + 1)x}\right)^{k-1} = \frac{w + 1}{w^k + 1} \geq 1.$$

Define

$$h_1(w, k) = \sqrt[k-1]{\frac{w + 1}{w^k + 1}} \geq 1.$$

Then, for every $k$ the derivative is zero at

$$x = \frac{d\, h_1(w, k)}{(w + 1)h_1(w, k) + 1} \leq 0.$$

For odd $k$ the derivative is zero also at

$$x = \frac{d\, h_1(w, k)}{(w + 1)h_1(w, k) - 1} \leq 0.$$

In the triangle $OBJ$ we have $x \geq y \geq 0$, so the derivative of $s_k$ has the same sign for $x \in \left[0, \frac{2d+1}{w+1}\right]$. The derivative at $x = 0$ equals $-k(w + 1)d^{k-1}$. Therefore, for odd (even) $k$ the derivative is negative (positive).

Thus, for $d \in \left[-\frac{1}{2}, 0\right]$ and $w \in [0,1]$ we proved the following: For odd (even) $k$ the maximum (minimum) value of $s_k$ over the line $y = wx$ is achieved at $x = 0$ and the minimum (maximum) value of $s_k$ over the line $y = wx$ is achieved at $x = \frac{2d+1}{w+1}$.

At $x = 0$ we have $s_k = 1 + d^k + (-d - 1)^k$ and at $x = \frac{2d+1}{w+1}$ we have $s_k = h_2(w, d, k)$, where

$$h_2(w, d, k) = 1 + (w^k + 1)\left(\frac{2d + 1}{w + 1}\right)^k + 2(-d - 1)^k.$$

The derivative of $h_2$ with respect to $w$ is

$$\frac{\partial h_2}{\partial w} = \left(-k\,\frac{w^k + 1}{(w + 1)^{k+1}} + k\,\frac{w^{k-1}}{(w + 1)^k}\right)(2d + 1)^k.$$

Since we already dealt with the case $d = -\frac{1}{2}$ we can assume $d \neq -\frac{1}{2}$. Therefore, the derivative is zero when

$$k\frac{w^{k-1}}{(w+1)^k} = k\ \frac{w^k+1}{(w+1)^{k+1}},$$

$$(w+1)w^{k-1} = \ w^k+1,$$

$$w^{k-1} = 1.$$

So, for every $k$ the derivative is zero at $w = 1$. For odd $k$ the derivative is zero also at $w = -1$. In any case, the derivative has the same sign in the range $(-1,1)$ of $w$.

Since we assume $-\frac{1}{2} < d \leq 0$ at $w = 0$ the derivative is negative and its value is $-k(2d+1)^k$. Thus, the minimum of $h_2$ is achieved at $w = 1$ and its value there is $1 + 2(d+\frac{1}{2})^k + 2(-d-1)^k$. The maximum of $h_2$ is achieved at $w = 0$ and its value there is $1 + (2d+1)^k + 2(-d-1)^k$.

This completes the proof of the lemma.

**Lemma 2**: Let $s_3 = 1 + x^3 + y^3 + (d-x-y)^3 + (-d-1)^3$. Then,

1. If $d \in \left[-\frac{1}{2}, -\frac{3}{4} + \frac{\sqrt{5}}{4}\right]$ then $s_3 \geq 0$ within the triangle $OBJ$,
2. If $d \in \left[-\frac{3}{4} + \frac{\sqrt{5}}{4}, 0\right]$ then $s_3 \geq 0$ within the 4-vertex shape $OHIJ$, which is formed by the intersection of the shape $P$ with the triangle $OBJ$.

**Proof**: Note that

$$\begin{aligned} s_3 &= 1 + x^3 + y^3 + (d-x-y)^3 + (-d-1)^3 \\ &= 3(d-x)y^2 + 3(2dx - d^2 - x^2)y + 3(dx^2 - d^2x - d - d^2). \end{aligned}$$

By Lemma 1 the minimum value over the triangle $OBJ$ is

$$\begin{aligned} 1 + 2\left(d+\frac{1}{2}\right)^3 + 2(-d-1)^3 &= -3d^2 - \frac{9}{2}d - \frac{3}{4} \\ &= -3\left(d - \left(-\frac{3}{4} + \frac{\sqrt{5}}{4}\right)\right)\left(d - \left(-\frac{3}{4} - \frac{\sqrt{5}}{4}\right)\right). \end{aligned}$$

Therefore, for $d \in \left[-\frac{1}{2}, -\frac{3}{4} + \frac{\sqrt{5}}{4}\right]$ we have $s_3 \geq 0$ over the entire triangle $OBJ$.

For $d \in \left[-\frac{3}{4} + \frac{\sqrt{5}}{4}, 0\right]$ we shall find when $s_3 \geq 0$.

When $d = 0$ we get $s_3 = -3xy(x+y)$. As $x \geq y \geq 0$ in the triangle we conclude that $s_3 \geq 0$ if and only if $y = 0$. We already know that the when $d = 0$ the shape $P$ becomes the triangle $OJG$, so the intersection with the triangle $OBJ$ is the line $OJ$. Therefore, the lemma is proved in this case.

Assume that $d < 0$. Then $x \geq 0 > d$ so $s_3$ is a quadratic function in the variable $y$.

The roots of $s_3$ as a function of $y$ are:

$$y_1 = -\frac{1}{2}(x-d) + \frac{1}{2}\sqrt{\frac{x^3 + dx^2 - d^2x - 4d - 4d^2 - d^3}{x-d}},$$

$$y_2 = -\frac{1}{2}(x-d) - \frac{1}{2}\sqrt{\frac{x^3 + dx^2 - d^2x - 4d - 4d^2 - d^3}{x-d}}.$$

Note that $y_1$ is exactly the function $h(t)$ defined at the beginning of this paper with $t$ replaced by $x$.

First we find when $y_1, y_2$ are real. As $x > d$ we require

$$h_3(x) = x^3 + dx^2 - d^2x - 4d - 4d^2 - d^3 \geq 0.$$

The derivative of $h_3$ with respect to $x$ is

$$\frac{\partial h_3}{\partial x} = 3x^2 + 2dx - d^2 = 3(x+d)\left(x - \frac{1}{3}d\right).$$

So $h_3$ has a local minimum at $x = -d$ and a local maximum at $x = \frac{1}{3}d$. Since $d < 0$ we get

$$h_3\left(\frac{1}{3}d\right) = -\frac{32}{27}d^3 - 4d^2 - 4d = -4d(\frac{8}{27}d^2 + d + 1) \geq 0,$$

$$h_3(-d) = -4d^2 - 4d = -4d(d+1) \geq 0.$$

Therefore, $h_3$ has a single real root $x_3(d) \leq \frac{1}{3}d < 0$. For $x \geq x_3(d)$ we have $h_3(x) \geq 0$.

By using the formula for the roots of a cubic equation we get

$$x_3(d) = -\frac{1}{3}d + \frac{2\sqrt[3]{4}d^2}{3\sqrt[3]{\Delta(d)}} + \frac{\sqrt[3]{2}\sqrt[3]{\Delta(d)}}{3},$$

where,

$$\Delta(d) = 4d^3 + 27d^2 + 27d + 3\sqrt{3}\sqrt{d^2(d+1)(8d^2 + 27d + 27)}.$$

Note that $\Delta(d) < 0$ for $-1 \leq d < 0$. Otherwise, $x_3(d) \geq -\frac{1}{3}d > 0 > \frac{1}{3}d \geq x_3(d)$, which is a contradiction. Also note that $\Delta(d)$ is exactly the function $r(d)$ defined at the beginning of this paper.

We conclude that $y_1, y_2$ are real for $x \geq x_3(d)$. In that case $y_1 \geq y_2$. As $x > d$ we have $y_2 \leq 0$ and when $y$ is in the range $[y_2, y_1]$ we have $s_3 \geq 0$.

Next we check for what values of $x$ the function $y_1$ intersects the triangle $OBJ$.

Notice that:

$$h_3(x) = (x-d)^3 + 4d(x+1)(x-d-1).$$

Inside the triangle $OBJ$ we have $0 \leq x \leq 2d+1$. As $d < 0$ we get $x - d - 1 \leq d < 0$. Therefore, $h_3(x) \geq (x-d)^3 \geq 0$. We conclude that $y_1 \geq 0$.

We first check when $y_1 \leq x$:

$$x \geq -\frac{1}{2}(x-d) + \frac{1}{2}\sqrt{\frac{x^3 + dx^2 - d^2x - 4d - 4d^2 - d^3}{x-d}},$$

$$\left(\frac{3}{2}x - \frac{1}{2}d\right)^2 \geq \frac{x^3 + dx^2 - d^2x - 4d - 4d^2 - d^3}{4(x-d)},$$

$$4\left(\frac{3}{2}x - \frac{1}{2}d\right)^2(x-d) - (x^3 + dx^2 - d^2x - 4d - 4d^2 - d^3) \geq 0,$$

$$8x^3 - 16dx^2 + 8d^2x + 4d + 4d^2 \geq 0.$$

Let

$$h_4(x) = 8x^3 - 16dx^2 + 8d^2x + 4d + 4d^2.$$

The derivative of $h_4$ with respect to $x$ is

$$\frac{\partial h_4}{\partial x} = 24x^2 - 32dx + 8d^2 = 24(x-d)\left(x - \frac{1}{3}d\right).$$

At $x = d$ there is a local maximum of $h_4$ and its value there is

$$h_4(d) = 4d + 4d^2 = 4d(d+1) < 0.$$

Therefore, $h_4$ has a single real root.

The derivative is positive for $x > \frac{1}{3}d$. Also, as $-\frac{3}{4} + \frac{\sqrt{5}}{4} \leq d < 0$,

$$h_4(0) = 4d + 4d^2 = 4d(d+1) < 0,$$

$$h_4\left(d+\frac{1}{2}\right)=4d^2+6d+1=4\left(d-\left(-\frac{3}{4}+\frac{\sqrt{5}}{4}\right)\right)\left(d-\left(-\frac{3}{4}-\frac{\sqrt{5}}{4}\right)\right)\geq 0.$$

We conclude that there is a single root of $h_4$ in the range $\left[0,d+\frac{1}{2}\right]$, which we denote by $x_4(d)$. By using the formula for the roots of a cubic equation we get

$$x_4(d)=\frac{2}{3}d-\frac{\sqrt[3]{4}d^2}{3\sqrt[3]{\Delta(d)}}-\frac{\sqrt[3]{\Delta(d)}}{3\sqrt[3]{4}},$$

where $\Delta(d)$ is as before. Also note that $x_4(d)$ is exactly the function $f(d)$ defined at the beginning of this paper.

We found that $y_1\leq x$ for $-\frac{3}{4}+\frac{\sqrt{5}}{4}\leq d<0$ and $x\in\left[x_4(d),d+\frac{1}{2}\right]$. In other words, for $x\in\left[f(d),d+\frac{1}{2}\right]$ the only pairs $(x,y)$ that lie in the triangle and meet the condition $s_3\geq 0$ are those that have $0\leq y\leq h(x)$. In particular, as we stated at the beginning of this paper:

$$h(f(d))=y_1(x_4(d))=x_4(d)=f(d).$$

Next we check when $y_1\leq -x+2d+1$:

$$-x+2d+1\geq -\frac{1}{2}(x-d)+\frac{1}{2}\sqrt{\frac{x^3+dx^2-d^2x-4d-4d^2-d^3}{x-d}},$$

$$\left(1+\frac{3}{2}d-\frac{1}{2}x\right)^2\geq\frac{x^3+dx^2-d^2x-4d-4d^2-d^3}{4(x-d)},$$

$$4\left(1+\frac{3}{2}d-\frac{1}{2}x\right)^2(x-d)-(x^3+dx^2-d^2x-4d-4d^2-d^3)\geq 0,$$

$$-4(2d+1)x^2+4(2\mathrm{d}+1)^2x-8d^2(d+1)\geq 0.$$

Let

$$h_5(x)=-4(2d+1)x^2+4(2\mathrm{d}+1)^2x-8d^2(d+1).$$

Note that $h_5$ is a quadratic since we assume $d\geq-\frac{3}{4}+\frac{\sqrt{5}}{4}>-\frac{1}{2}$. The roots of $h_5$ are:

$$p_1(d)=d+\frac{1}{2}+\sqrt{\frac{\left(d-\left(-\frac{3}{4}+\frac{\sqrt{5}}{4}\right)\right)\left(d-\left(-\frac{3}{4}-\frac{\sqrt{5}}{4}\right)\right)}{2d+1}},$$

$$p_2(d) = d + \frac{1}{2} - \sqrt{\frac{\left(d - \left(-\frac{3}{4} + \frac{\sqrt{5}}{4}\right)\right)\left(d - \left(-\frac{3}{4} - \frac{\sqrt{5}}{4}\right)\right)}{2d+1}}.$$

Note that $p_1(d) = d + \frac{1}{2} + g(d)$, where $g(d)$ is the function defined at the beginning of this paper.

These roots are real since we deal with the case $d \geq -\frac{3}{4} + \frac{\sqrt{5}}{4}$. Also, $p_1(d) \geq d + \frac{1}{2} \geq p_2(d)$, and $h_5(x) \geq 0$ if and only if $x \in [p_2(d), p_1(d)]$.

Finally, since $h_5(2d+1) = -8d^2(d+1) < 0$ and $2d + 1 \geq d + \frac{1}{2}$ we get that $p_1(d) \leq 2d + 1$.

We found that $y_1 \leq -x + 2d + 1$ for $-\frac{3}{4} + \frac{\sqrt{5}}{4} \leq d < 0$ and $x \in \left[d + \frac{1}{2}, p_1(d)\right]$. In other words, for $x \in \left[d + \frac{1}{2}, d + \frac{1}{2} + g(d)\right]$ the only pairs $(x, y)$ that lie in the triangle and meet the condition $s_3 \geq 0$ are those that have $0 \leq y \leq h(x)$. In particular, as we stated at the beginning of this paper:

$$h\left(d + \frac{1}{2} + g(d)\right) = y_1(p_1(d)) = -p_1(d) + 2d + 1 = d + \frac{1}{2} - g(d).$$

This completes the proof of the lemma.

The final result we shall need is:

**Lemma 3**: Define the following symmetric $5 \times 5$ matrix families for $d \in \left[-\frac{1}{2}, 0\right]$:

$$A(x) = \begin{pmatrix} 0 & 0 & f_1 & 0 & f_1 \\ 0 & 0 & 0 & 0 & g_1 \\ f_1 & 0 & 0 & g_1 & -d \\ 0 & 0 & g_1 & 0 & 0 \\ f_1 & g_1 & -d & 0 & 0 \end{pmatrix},$$

where,

$$f_1 = f_1(x, d) = \sqrt{\frac{1}{2}(x+1)(d+1-x)},$$

$$g_1 = g_1(x, d) = \sqrt{x(x-d)}.$$

$$B(x,y)=\frac{1}{2u^2}\begin{pmatrix} 0 & 0 & 2u^3 & 0 & 2u^3 \\ 0 & 0 & 0 & (d+1)y & v \\ 2u^3 & 0 & 0 & v & \frac{1}{3}s_3 \\ 0 & (d+1)y & v & 0 & 0 \\ 2u^3 & v & \frac{1}{3}s_3 & 0 & 0 \end{pmatrix},$$

where,

$$s_3=s_3(x,y,d)=1+x^3+y^3+(d-x-y)^3+(-d-1)^3,$$

$$u=u(x,y,d)=\sqrt{\frac{(w_1-x)(x-w_2)(x+y)}{2x}},$$

$$v=v(x,y,d)$$
$$=\sqrt{(x+y)(x+1)(x-d)(x-d-1)(x+y+1)(x+y-d-1)},$$

$$w_1=w_1(x,y,d)=\frac{1}{2}d+\frac{1}{2}\sqrt{d^2+\frac{4(d+1)x}{x+y}},$$

$$w_2=w_2(x,y,d)=\frac{1}{2}d-\frac{1}{2}\sqrt{d^2+\frac{4(d+1)x}{x+y}}.$$

Then

1. $A(x)$ is nonnegative for any $x$ within the line segment $OJ$,
2. If $d\in\left[-\frac{1}{2},-\frac{3}{4}+\frac{\sqrt{5}}{4}\right]$ then $B(x,y)$ is nonnegative for any pair $(x,y)$ with $y>0$ within the triangle $OBJ$,
3. If $d\in\left[-\frac{3}{4}+\frac{\sqrt{5}}{4},0\right]$ then $B(x,y)$ is nonnegative for any pair $(x,y)$ with $y>0$ within the 4-vertex shape $OHIJ$, which is formed by the intersection of the shape $P$ with the triangle $OBJ$,
4. The spectrum of $A(x)$ is $\sigma=(1,x,0,d-x,-d-1)$ for any $x$ within the line segment $OJ$,
5. The spectrum of $B(x,y)$ is $\sigma=(1,x,y,d-x-y,-d-1)$ for any pair $(x,y)$ with $y>0$ within the triangle $OBJ$.

**Proof**: A point $(x,y)$ within the triangle $OBJ$ satisfies the following conditions:

$$y\le x,$$

$$y\le -x+2d+1,$$

$$y\ge 0.$$

Using these conditions and the fact that $d \leq 0$ it is trivial to show that all the matrix elements of $A(x)$ are nonnegative.

For $y > 0$ we get from the above conditions $0 < x < 2d + 1$. Also, as $y > 0$ we have $d \neq -\frac{1}{2}$. Since $-\frac{1}{2} < d \leq 0$, we get

$$d + 1 + d(2d + 1) = (2d + 1)^2 - 2d(d + 1) \geq (2d + 1)^2 > (2d + 1)x.$$

Therefore,

$$\frac{d + 1}{x - d} > 2d + 1.$$

Using this result and by the second triangle condition

$$y \leq -x + 2d + 1 < -x + \frac{d + 1}{x - d},$$

$$x + y < \frac{d + 1}{x - d},$$

$$4x(x - d) < 4x\left(\frac{d + 1}{x + y}\right),$$

$$(2x - d)^2 < d^2 + \frac{4(d + 1)x}{x + y},$$

$$0 < 2x - d < \sqrt{d^2 + \frac{4(d + 1)x}{x + y}},$$

$$x < \frac{1}{2}d + \frac{1}{2}\sqrt{d^2 + \frac{4(d + 1)x}{x + y}} = w_1.$$

Therefore, we showed that $x < w_1$ for a point $(x, y)$ with $y > 0$ within the triangle $OBJ$. By the triangle conditions, $y > 0$, and $-\frac{1}{2} < d \leq 0$ it can be easily shown that $u > 0$ and that all the matrix elements of $B(x, y)$, apart from $\frac{1}{3}s_3$, are nonnegative. By Lemma 2 we know that $s_3 \geq 0$ within the given regions for which $B(x, y)$ is claimed to be nonnegative.

The characteristic polynomial of $A(x)$ is:

$$
\begin{aligned}
p_A(z) &= z^5 - \left(2{g_1}^2 + 2{f_1}^2 + d^2\right)z^3 + 2d{f_1}^2z^2 + \left({g_1}^4 + 2{f_1}^2{g_1}^2\right)z \\
&= z^5 - (2x(x-d) + (x+1)(d+1-x) + d^2)z^3 \\
&+ d(x+1)(d+1-x)z^2 \\
&+ \left(x^2(x-d)^2 + x(x-d)(x+1)(d+1-x)\right)z \\
&= z^5 - (x^2 - dx + d^2 + d + 1)z^3 + (-dx^2 + d^2x + d^2 + d)z^2 \\
&+ (x^2 + dx^2 - d^2x - dx)z \\
&= (z-1)(z-x)z\big(z-(d-x)\big)\big(z-(-d-1)\big).
\end{aligned}
$$

This shows that the spectrum of $A(x)$ is $\sigma = (1, x, 0, d-x, -d-1)$.

We shall now compare the coefficients of the characteristic polynomial of $B(x,y)$ and the coefficients of the polynomial

$$
\begin{aligned}
q(z) &= (z-1)(z-x)(z-y)\big(z-(d-x-y)\big)\big(z-(-d-1)\big) \\
&= z^5 + q_4z^4 + q_3z^3 + q_2z^2 + q_1z + q_0.
\end{aligned}
$$

As before let

$$s_k = 1 + x^k + y^k + (d-x-y)^k + (-d-1)^k.$$

Obviously, $s_1 = 0$ and $q_4 = 0$.

By the Newton identities:

$$2q_3 = -q_4s_1 - s_2 = -s_2,$$

$$-3q_2 = q_3s_1 - q_4s_2 + s_3 = s_3,$$

$$4q_1 = -q_2s_1 - q_3s_2 - q_4s_3 - s_4 = -q_3s_2 - s_4,$$

$$-5q_0 = q_1s_1 + q_2s_2 + q_3s_3 + q_4s_4 + s_5 = q_2s_2 + q_3s_3 + s_5.$$

We therefore get:

$$q_3 = -\frac{1}{2}s_2,$$

$$q_2 = -\frac{1}{3}s_3,$$

$$q_1 = -\frac{1}{4}s_4 + \frac{1}{8}{s_2}^2,$$

$$q_0 = \frac{1}{6}s_2s_3 - \frac{1}{5}s_5.$$

The characteristic polynomial of $B(x,y)$ is:

$$p_B(z) = z^5 + p_3z^3 + p_2z^2 + p_1z + p_0,$$

where,

$$p_3 = -\left(\frac{(d+1)y}{2u^2}\right)^2 - 2\left(\frac{v}{2u^2}\right)^2 - 2u^2 - \left(\frac{s_3}{6u^2}\right)^2$$
$$= -\frac{1}{4u^4}\left((d+1)^2y^2 + 2v^2 + 8u^6 + \frac{1}{9}{s_3}^2\right),$$

$$p_2 = -2u^2\frac{s_3}{6u^2} = -\frac{1}{3}s_3 = q_2,$$

$$p_1 = -2\frac{(d+1)y}{2u^2}\left(\frac{v}{2u^2}\right)^2\frac{s_3}{6u^2} + \left(\frac{v}{2u^2}\right)^4 + 2\left(\frac{(d+1)y}{2u^2}\right)^2 u^2$$
$$+ \left(\frac{(d+1)y}{2u^2}\right)^2\left(\frac{s_3}{6u^2}\right)^2 + 2\left(\frac{v}{2u^2}\right)^2 u^2$$
$$= \frac{1}{16u^8}\left(-\frac{2}{3}(d+1)yv^2s_3 + v^4 + 8(d+1)^2y^2u^6\right.$$
$$\left.+\frac{1}{9}(d+1)^2y^2{s_3}^2 + 8v^2u^6\right),$$

$$p_0 = -2\frac{(d+1)y}{2u^2}\left(\frac{v}{2u^2}\right)^2 u^2 + 2\left(\frac{(d+1)y}{2u^2}\right)^2 u^2\frac{s_3}{6u^2}$$
$$= \frac{1}{4u^4}\left(-(d+1)yv^2 + \frac{1}{3}(d+1)^2y^2s_3\right).$$

Substituting $w_1$ and $w_2$ in the function $u$ gives:

$$u = \sqrt{\frac{1}{2}\left((d-x)(x+y)+d+1\right)} = \sqrt{\frac{1}{2}(-x^2 - xy + dx + dy + d + 1)}.$$

The following expressions are polynomials in $x$ and $y$ since $u$ and $v$ appear only with even powers. Therefore, it is straightforward to check that

$$-4u^4p_3 = (d+1)^2y^2 + 2v^2 + 8u^6 + \frac{1}{9}{s_3}^2 = 2u^4s_2 = -4u^4q_3,$$

$$16u^8p_1 = -\frac{2}{3}(d+1)yv^2s_3 + v^4 + 8(d+1)^2y^2u^6 + \frac{1}{9}(d+1)^2y^2{s_3}^2$$
$$+ 8v^2u^6 = u^8(-4s_4 + 2{s_2}^2) = 16u^8q_1,$$

$$4u^4p_0 = -(d+1)yv^2 + \frac{1}{3}(d+1)^2y^2s_3 = 4u^4\left(\frac{1}{6}s_2s_3 - \frac{1}{5}s_5\right) = 4u^4q_0.$$

Since $u > 0$ we get that all the coefficient of $q(z)$ and $p_B(z)$ are equal, so the spectrum of $B(x,y)$ is $\sigma = (1, x, y, d-x-y, -d-1)$.

This completes the proof of the lemma.

## 5. Proof of main result

**Proof of Theorem 2 (Main result)**: We shall prove the theorem by considering three different cases:

1. $x \leq 0$,
2. $x > 0$ and $y \leq 0$,
3. $x > 0$ and $y > 0$.

Note that the triangle $ABC$ for $d \in \left[-\frac{3}{4}, -\frac{1}{2}\right]$ is fully covered by case 1 and case 2.

For $d \in \left[-\frac{1}{2}, -\frac{3}{4} + \frac{\sqrt{5}}{4}\right]$ case 3 is the triangle $OBJ$ without the edge $OJ$. For $d \in \left[-\frac{3}{4} + \frac{\sqrt{5}}{4}, 0\right]$ case 3 is the 4-vertex shape $OHIJ$, which is formed by the intersection of the shape $P$ with the triangle $OBJ$, without the edge $OJ$. Note that when $d = 0$ this intersection is empty so case 3 is never happens.

Let $d \in \left[-\frac{3}{4}, -\frac{1}{2}\right]$. We already know that if $\sigma = (1, x, y, d - x - y, -d - 1)$ is a normalized spectrum of a matrix $A \in \widetilde{\mathfrak{R}}$ then $(x, y)$ must lie within the triangle $ABC$. This proves the necessity of the condition.

To prove sufficiency, assume that $(x, y)$ lies within the triangle $ABC$. As mentioned before, case 3 is impossible, so we are left with the other two cases.

The triangle $ABC$ conditions are:

$$y \leq x,$$

$$y \leq -x + 2d + 1,$$

$$y \geq \frac{1}{2}(d - x).$$

If $x \leq 0$ then the conditions of Theorem 4 are satisfied and therefore $\sigma = (1, x, y, d - x - y, -d - 1)$ is realized by a symmetric nonnegative $5 \times 5$ matrix $A$. Since the sum of the elements of $\sigma$ are zero, then $A \in \widetilde{\mathfrak{R}}$.

Note that this proof is valid for $d \in \left[-\frac{3}{4}, 0\right]$.

If $x > 0$ and $y \leq 0$ then let $\sigma = (1, x, y, d - x - y, -d - 1) = (\lambda_1, \lambda_2, \lambda_3, \lambda_4, \lambda_5)$ and $K_1 = \{3,5\}, K_2 = \{4\}$. We have $\sum_{i=1}^{5} \lambda_i = 0$, $\sum_{i \in K_1} \lambda_i = y - d - 1$ and $\sum_{i \in K_2} \lambda_i = d - x - y$.

Assume that $-\sum_{i\in K_1}\lambda_i < -\sum_{i\in K_2}\lambda_i$. In that case $d+1-y < x+y-d$. Therefore, $2y > -x+2d+1$. As $y \le 0$ we have $2y \le y$. Using the second triangle condition we get $-x+2d+1 \ge y \ge 2y > -x+2d+1$, which is a contradiction and therefore, $-\sum_{i\in K_1}\lambda_i \ge -\sum_{i\in K_2}\lambda_i$.

As we assumed that $d \le -\frac{1}{2}$ we have $3d+2 \le -d$. From the second and third triangle conditions we have $x \le 3d+2$ and $y \ge \frac{1}{2}(d-x)$. Therefore, $y \ge \frac{1}{2}(d-x) \ge \frac{1}{2}d - \frac{1}{2}(3d+2) \ge \frac{1}{2}d + \frac{1}{2}d = d$. This means that $1 \ge d+1-y = -\sum_{i\in K_1}\lambda_i$.

All the conditions of Theorem 5 are satisfied and therefore we conclude that $\sigma = (1,x,y,d-x-y,-d-1)$ is realized by symmetric nonnegative $5\times 5$ matrix $A$. Since the sum of the elements of $\sigma$ are zero, then $A \in \widetilde{\mathfrak{R}}$.

Let $d \in \left[-\frac{1}{2}, -\frac{3}{4}+\frac{\sqrt{5}}{4}\right]$. We already know that if $\sigma = (1,x,y,d-x-y,-d-1)$ is a normalized spectrum of a matrix $A \in \widetilde{\mathfrak{R}}$ then $(x,y)$ must lie within the quadrangle $ABFG$. This proves the necessity of the condition.

To prove sufficiency, assume that $(x,y)$ lies within the quadrangle $ABFG$. We deal with case 1 and case 2 and leave case 3 for later.

The quadrangle $ABFG$ conditions are:

$$y \le x,$$

$$y \le -x+2d+1,$$

$$y \ge \frac{1}{2}(d-x),$$

$$x \le d+1.$$

Obviously these are triangle $ABC$ conditions plus the MN condition.

Case 1 is already proved as noted above.

Case 2 is proved using Theorem 5 as before. We assume $x > 0$ and $y \le 0$.

If $x > 0$ and $y \le 0$ then let $\sigma = (1,x,y,d-x-y,-d-1) = (\lambda_1,\lambda_2,\lambda_3,\lambda_4,\lambda_5)$ and $K_1 = \{3,4\}, K_2 = \{5\}$. We have $\sum_{i=1}^{5}\lambda_i = 0, \sum_{i\in K_1}\lambda_i = d-x, \sum_{i\in K_2}\lambda_i = -d-1$.

First assume that $x > 2d+1$. As $d \ge -\frac{1}{2}$ we get $x > 0$.

By the fourth quadrangle inequality $x \le d+1$, so $1 \ge x-d = -\sum_{i\in K_1}\lambda_i$.

By $x > 2d + 1$ we get $-\sum_{i \in K_1} \lambda_i = x - d > d + 1 = -\sum_{i \in K_2} \lambda_i$.

All the conditions of Theorem 5 are satisfied and therefore we conclude that $\sigma = (1, x, y, d - x - y, -d - 1)$ is realized by symmetric nonnegative $5 \times 5$ matrix $A$. Since the sum of the elements of $\sigma$ are zero, then $A \in \widetilde{\mathfrak{R}}$.

Next assume that $0 < x \leq 2d + 1$.

As $d \leq 0$ we have $1 \geq d + 1 = -\sum_{i \in K_2} \lambda_i$.

By $x \leq 2d + 1$ we get $-\sum_{i \in K_1} \lambda_i = x - d \leq d + 1 = -\sum_{i \in K_2} \lambda_i$.

Again, all the conditions of Theorem 5 are satisfied (with the roles of $K_1, K_2$ switched) and therefore $\sigma = (1, x, y, d - x - y, -d - 1)$ is realized by symmetric nonnegative $5 \times 5$ matrix $A$. Since the sum of the elements of $\sigma$ are zero, then $A \in \widetilde{\mathfrak{R}}$.

Note that this proof is valid for $d \in \left[-\frac{1}{2}, 0\right]$.

Let $d \in \left[-\frac{3}{4} + \frac{\sqrt{5}}{4}, 0\right]$. By Lemma 2 we know that if $\sigma = (1, x, y, d - x - y, -d - 1)$ is a normalized spectrum of a matrix $A \in \widetilde{\mathfrak{R}}$ then $(x, y)$ must lie within the shape $P$. This proves the necessity of the condition.

To prove sufficiency, assume that $(x, y)$ lies within the shape $P$. Case 1 and case 2 are already proved for this range of $d$ as noted above.

By Lemma 3 we know that for any pair $(x, y)$ which meets case 3 for $d \in \left[-\frac{1}{2}, 0\right]$ there is a matrix $B(x, y) \in \widetilde{\mathfrak{R}}$ with a spectrum $\sigma = (1, x, y, d - x - y, -d - 1)$. Therefore, we proved sufficiency.

This completes the proof of the theorem.

## Acknowledgement

I wish to thank Raphi Loewy for his helpful comments after reading an earlier draft of this paper.

## Appendix A

The following figures show the points defined in Section 2 for different values of $d$.

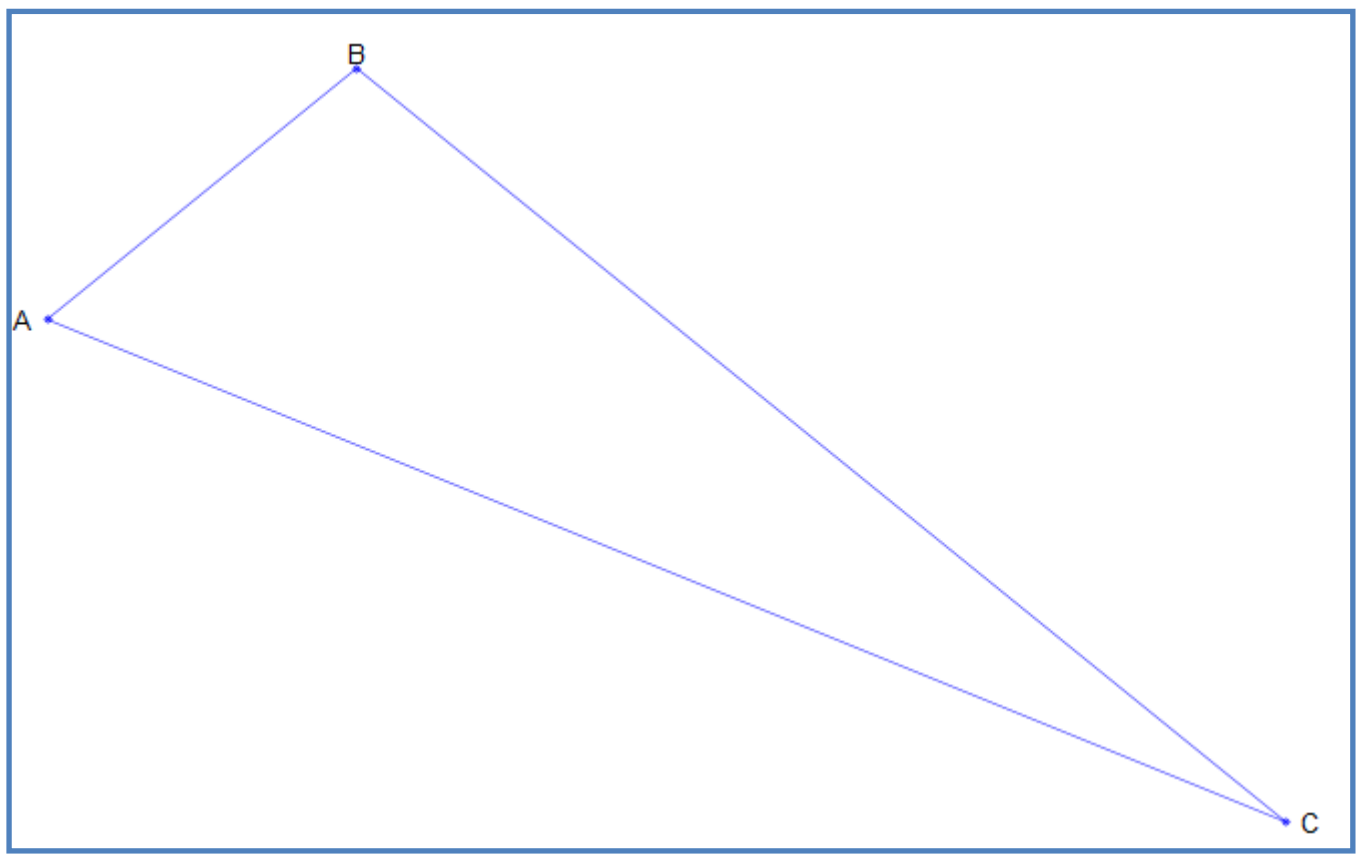

**Figure 1 –** $d \in \left(-\frac{3}{4}, -\frac{1}{2}\right)$

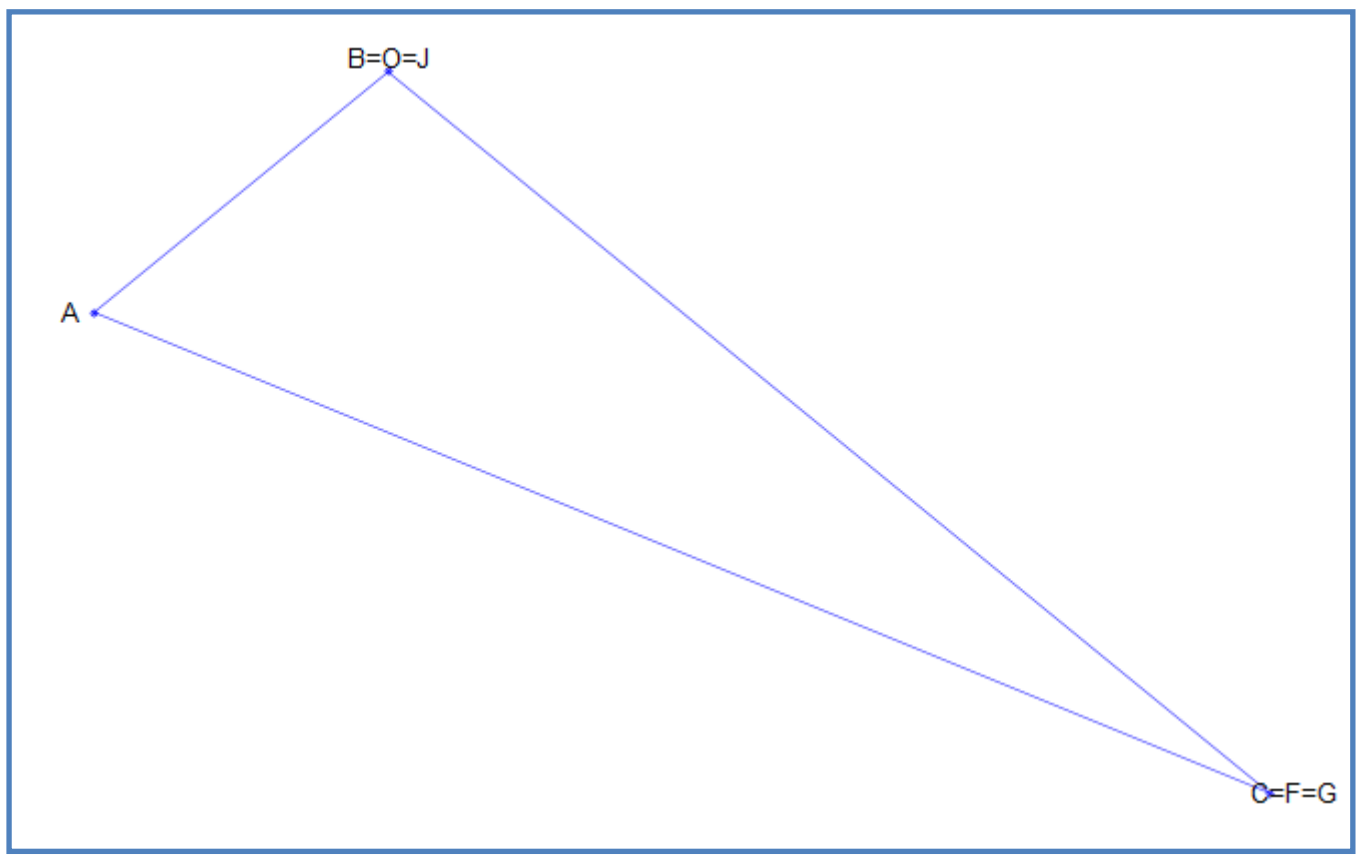

**Figure 2 –** $d = -\frac{1}{2}$

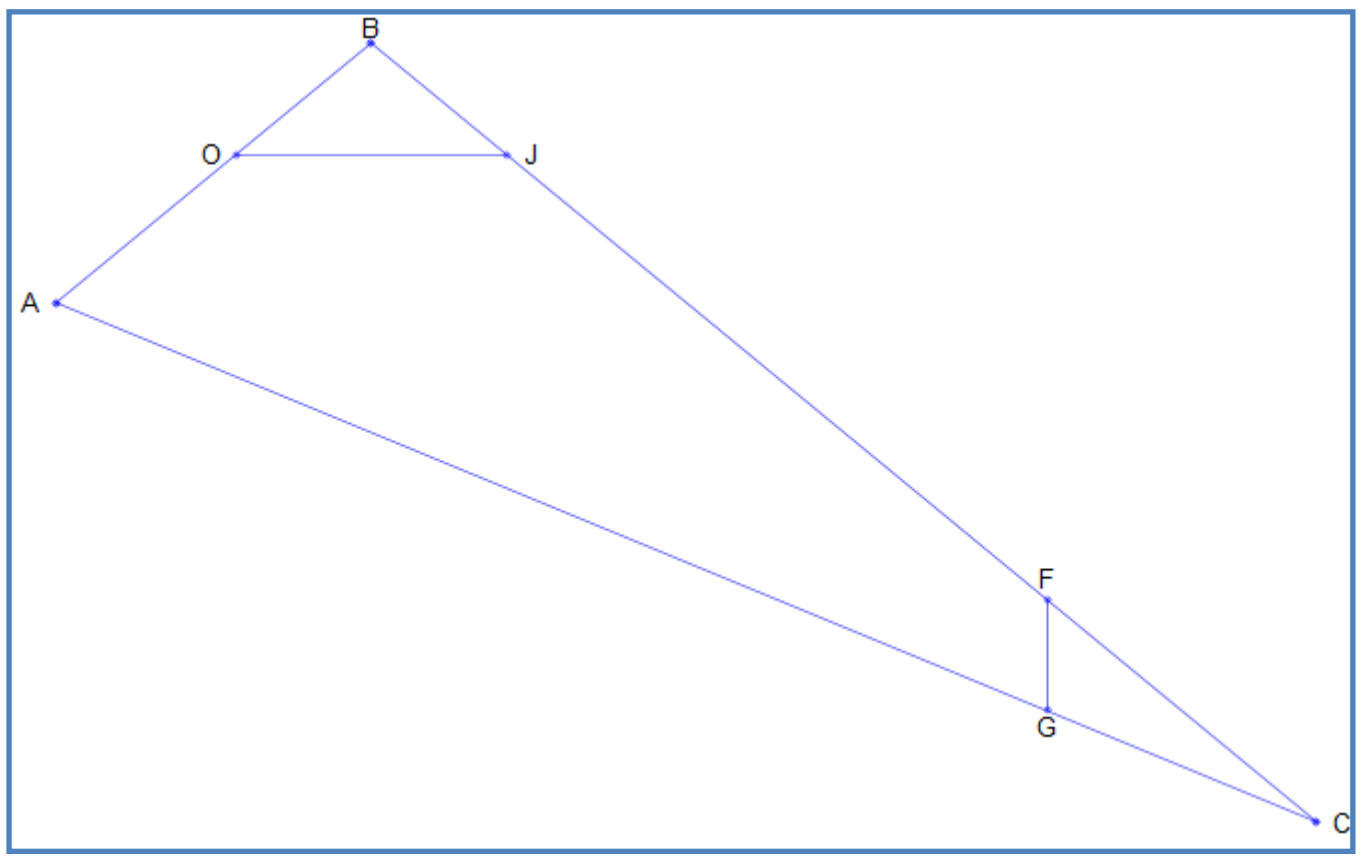

**Figure 3 –** $d \in \left(-\frac{1}{2}, -\frac{1}{3}\right)$

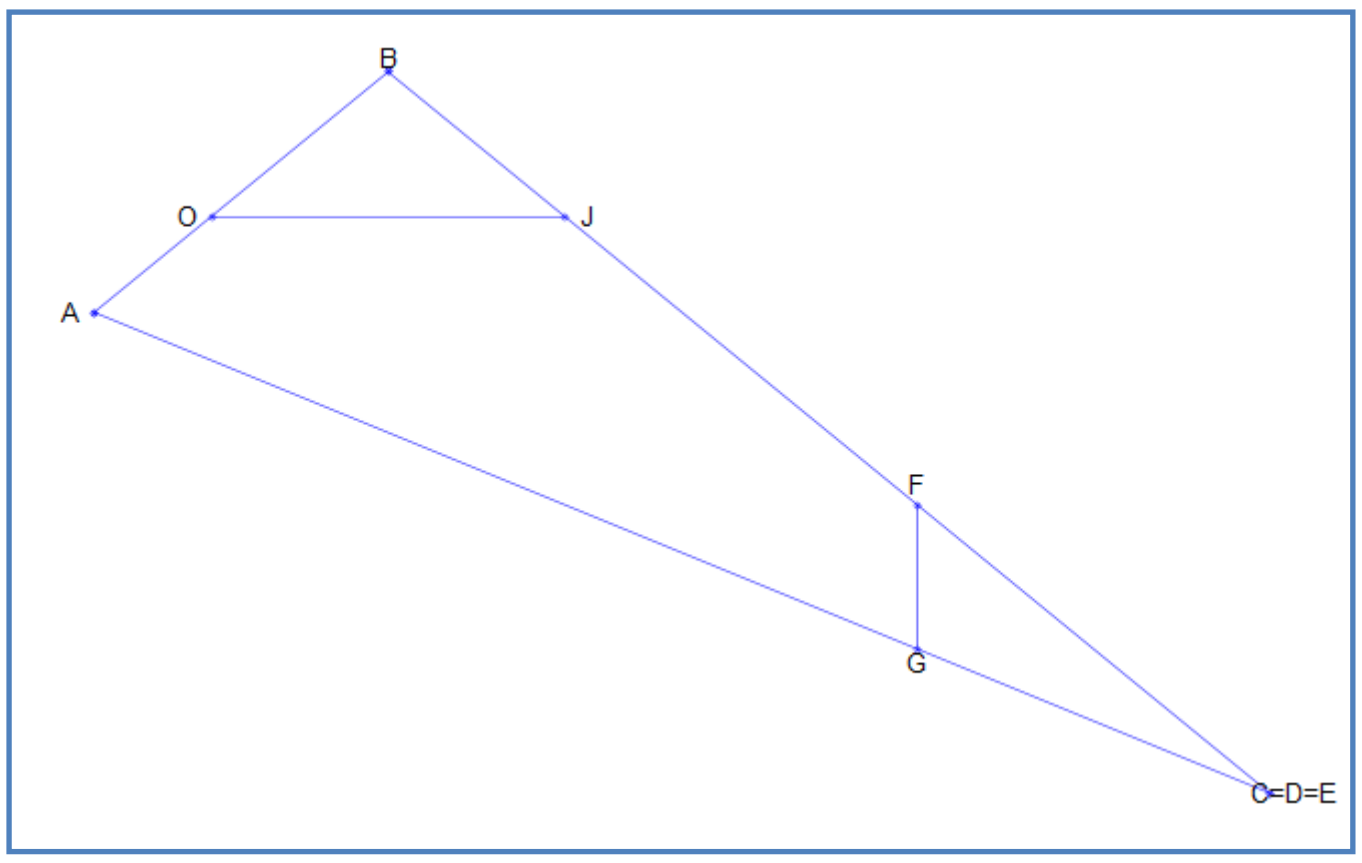

Figure 4 – $d = -\frac{1}{3}$

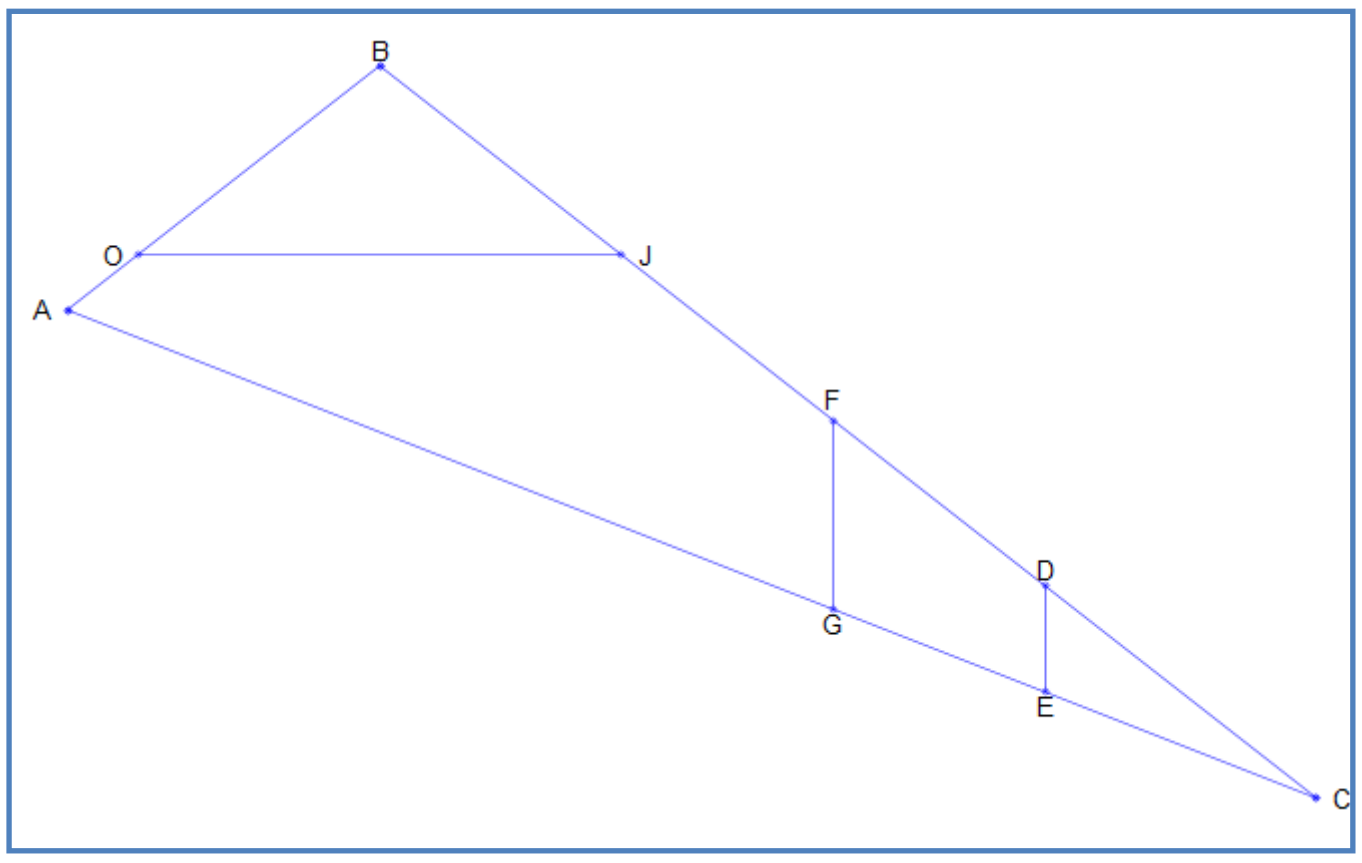

Figure 5 – $d \in \left(-\frac{1}{3}, -\frac{3}{4} + \frac{\sqrt{5}}{4}\right)$

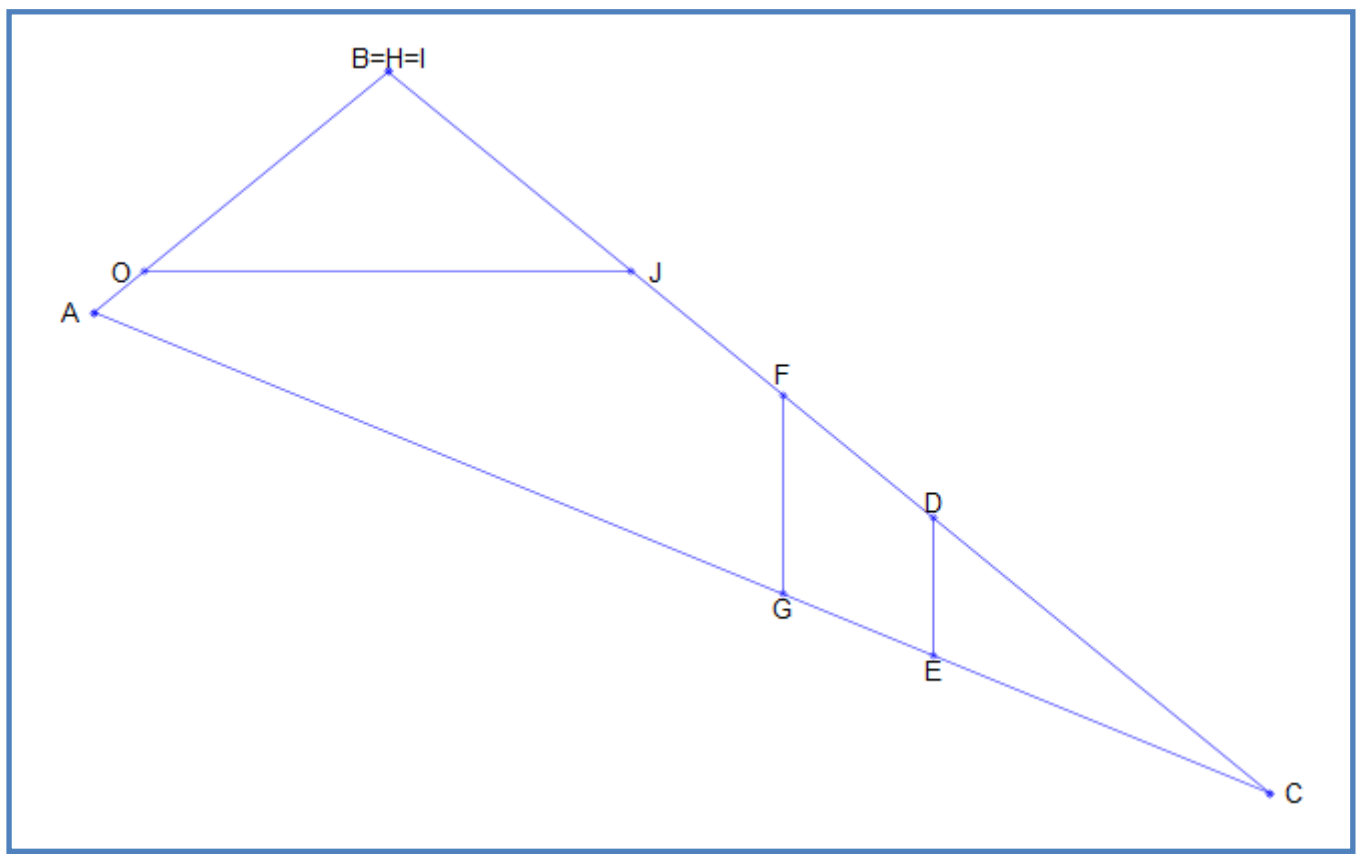

Figure 6 – $d = -\frac{3}{4} + \frac{\sqrt{5}}{4}$

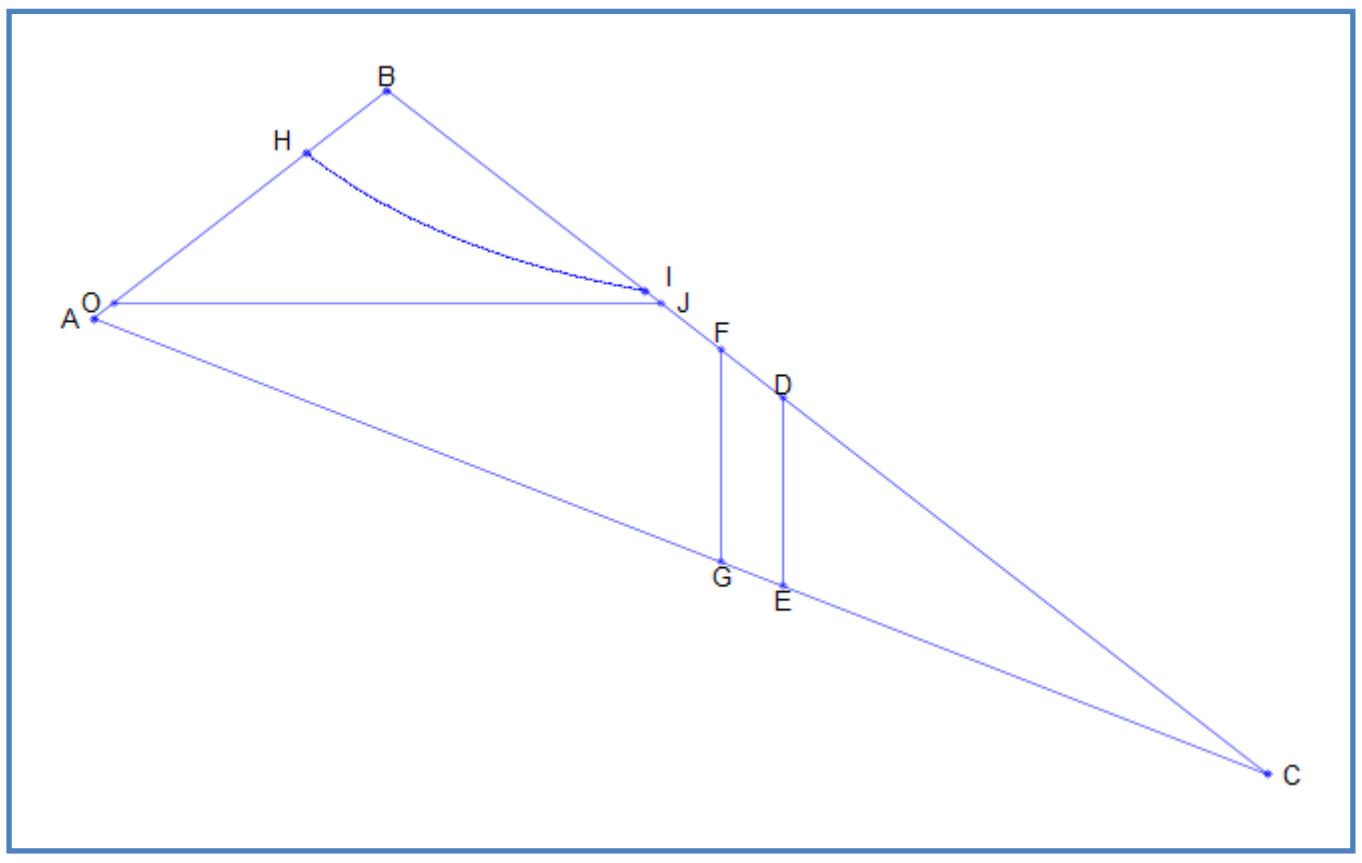

Figure 7 – $d \in \left(-\frac{3}{4}+\frac{\sqrt{5}}{4}, 0\right)$

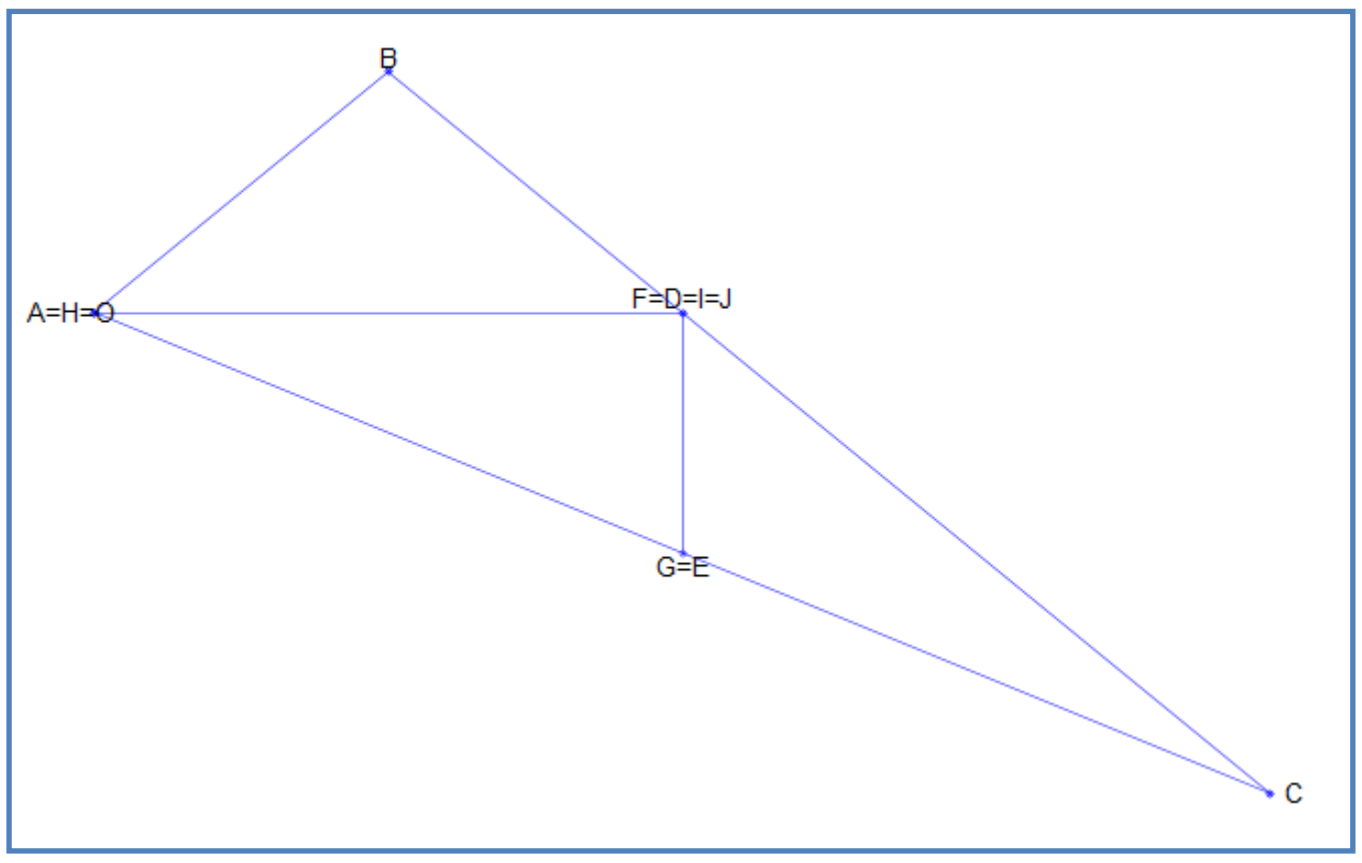

Figure 8 – $d = 0$